\documentclass{amsart}
\usepackage{amsmath,amssymb,amscd,amsthm,amsxtra,array}
\usepackage[
colorlinks=true,
linkcolor=blue,
 citecolor=blue,
  urlcolor=blue,
]{hyperref}

 \newtheorem{thm}{Theorem}[section]
 
 \newtheorem{lem}[thm]{Lemma}
 
 \theoremstyle{definition}
 
 \theoremstyle{remark}

 \numberwithin{equation}{section}

\numberwithin{equation}{section}

\begin{document}

\title[Improving bounds for singular operators via RHI for $A_{\infty}$]{Improving bounds for singular operators via Sharp Reverse H\"older Inequality for $A_{\infty}$
}

\author{Carmen Ortiz-Caraballo}
\address{Departamento de Matem\'aticas
\\ Escuela Polit\'ecnica, Universidad de Extremadura\\ Avda. Universidad, s/n, 10003 C\'aceres, Spain} \email{carortiz@unex.es}

\author{Carlos P\'erez}
\address{Departamento de An\'alisis Matem\'atico,
Facultad de Matem\'aticas, Universidad de Sevilla, 41080 Sevilla,
Spain} \email{carlosperez@us.es}

\author{Ezequiel Rela}
\address{Departamento de An\'alisis Matem\'atico,
Facultad de Matem\'aticas, Universidad de Sevilla, 41080 Sevilla,
Spain} \email{erela@us.es}
\thanks{The second author is supported by the Spanish Ministry of Science and Innovation grant MTM2009-08934,
the second and third authors are also supported by the Junta de Andaluc\'ia, grant FQM-4745.}



\date{\today}
\dedicatory{ Dedicated To Prof. Samko
}

\begin{abstract}
In this expository article we collect and discuss some recent results on different consequences of a Sharp Reverse H\"older Inequality for $A_{\infty}$ weights. For two given operators $T$ and $S$, we study $L^p(w)$ bounds of Coifman-Fefferman type:
$$ \|Tf\|_{L^p(w)}\le c_{n,w,p} \|Sf\|_{L^p(w)}, $$
that can be understood as a way to control $T$ by $S$.

We will focus on a \emph{quantitative} analysis of the constants involved and show that we can improve classical results regarding the dependence on the weight $w$ in terms of Wilson's $A_{\infty}$ constant
\begin{equation*}
  [w]_{A_{\infty}}:=\sup_Q\frac{1}{w(Q)}\int_Q M(w\chi_Q).
\end{equation*}

We will also exhibit recent improvements on the problem of finding sharp constants for weighted norm inequalities involving several singular opera\-tors. In the same spirit as in \cite{HP}, we obtain mixed $A_{1}$--$A_{\infty}$ estimates for the commutator $[b,T]$ and for its higher order analogue $T^k_{b}$.
A common ingredient in the proofs presented here is a recent improvement of the Reverse H\"older Inequality for $A_{\infty}$ weights involving Wilson's constant from \cite{HP}.

\end{abstract}

\maketitle

\section{Introduction}

The purpose of this survey is to exhibit and discuss some recent results involving improvements in two closely related topics in weight theory:  Coifman--Fefferman type inequalities for weighted $L^p$ spaces and boundedness of Calder\'on--Zygmund (C--Z) operators and  their commutators with BMO functions. More generally, we present several examples of the so called ``Calder\'on--Zygmund Principle'' which says, essentially, that any C--Z operator is controlled (in some sense) by an adequate maximal operator.  Typically, those inequalities provide  control over an operator $T$ with some singularity, like a singular integral operator, by means of a  ``better'' operator $S$, like a maximal operator.  As a model example of this phenomenom, we can take the classical Coifman--Fefferman  inequality involving a C--Z operator and the usual Hardy-Littlewood maximal function $M$ (see \cite{C}, \cite{CF}). Throughout this paper, we will denote  $L^p(\mathbb{R}^n,w)$ as $L^p(w)$. As usual, for $1<p<\infty$, we will denote  with $p'$ the dual index of $p$, defined by the equation $\frac{1}{p}+\frac{1}{p'}=1$. 

\begin{thm}[Coifman-Fefferman]\label{CF}
For any weight $w$ in the Muckenhoupt class $A_{\infty}$, the following norm inequality holds:
\begin{equation} \label{roughC-F}
\|Tf\|_{L^p(w)}\le c \,\|Mf\|_{L^p(w)},
\end{equation}
where $0<p<\infty$ and $c=c_{n,w,p}$ is a positive constant depending on the dimension $n$, the exponent $p$ and on the weight $w$.
\end{thm}

Let us recall that a weight $w$ is a non-negative measurable function. Such a function $w$ belongs to the  Muckenhoupt class $A_{p}$, $1<p<\infty$ if
\begin{equation*}
[w]_{A_p}:=\sup_{Q}\left(\frac{1}{|Q|}\int_{Q}w(y)\ dy \right)\left(\frac{1}{|Q|}\int_{Q}w(y)^{1-p'}\ dy \right)^{p-1}<\infty,
\end{equation*}
where the supremum is taken over all the cubes in $\mathbb{R}^n$. This number is called the $A_{p}$ constant or characteristic of the weight $w$. For $p=1$, the condition is that there exists a constant $c>0$ such that the Hardy-Littlewood maximal function $M$ satisfies the bound
\begin{equation*}
Mw(x)\le c \,w(x) \ \text{a.e.} \ x\in\mathbb{R}^n.
\end{equation*}
In that case, $[w]_{A_{1}}$ will denote the smallest of these constants. Since the $A_{p}$ classes are increasing with respect to $p$, we can define the $A_{\infty}$ class in the natural way by
 $$A_{\infty}:=\bigcup_{p>1}A_p.$$
For any weight in this larger class, the $A_{\infty}$ constant can be defined as follows:
\begin{equation*}
 \|w\|_{A_\infty}:=\sup_Q\Big(\frac{1}{|Q|}\int_Q w\Big)\exp\Big(\frac{1}{|Q|}\int_Q \log w^{-1}\Big).
\end{equation*}

This constant was introduced by Hru\v{s}\v{c}ev \cite{Hruscev} (see also \cite{GCRdF}) and   has been the standard $A_{\infty}$ constant until very recently when a  ``new'' $A_{\infty}$ constant was found to be better suited. This new constant is defined as
\begin{equation*}
  [w]_{A_\infty}:=\sup_Q\frac{1}{w(Q)}\int_Q M(w\chi_Q).
\end{equation*}
However, this constant was introduced  by M. Wilson a long time ago (see \cite{Wilson:87,Wilson:89,  Wilson-LNM}) with a different notation. This constant  is more relevant since there are examples of weights $w \in A_{\infty}$ so that   $[w]_{A_\infty} $ is much smaller than
$\|w\|_{A_\infty}.$ Indeed, first it can be shown that
\begin{equation}\label{mono}
  c_n\, [w]_{A_\infty}\leq \|w\|_{A_\infty}\leq[w]_{A_p},    \quad 1<p<\infty,
\end{equation}
where $c_n$ is a constant depending only on the dimension. The first inequality is the only nontrivial part and can be found in \cite{HP} where a more interesting fact is  also shown , namely that this inequality can be strict. More precisely, the authors exhibit a  family of weights $\{w_{t}\}$ such that   $[w_{t}]_{A_\infty}\le 4\log(t)$ and $\|w_{t}\|_{A_\infty}\sim t/\log(t)$ for $t\gg 1$.

It should be mentioned that this constant has also been used by Lerner \cite[Section~5.5]{L2} (see also \cite{L}) where the term ``$A_{\infty}$  constant'' was coined.

Going back to Theorem \ref{CF}, we can see that the result is, in the way that it is stated, a \emph{qualitative} result. It says that the maximal operator   $M$ acts as a ``control operator'' for C--Z operators, but the dependence of the constant $c$ on both $w$ and $p$ is not precise enough for some applications.   For instance, let us mention that a precise knowledge of the behavior of this constant was crucial in the proof of the following theorems from \cite{LOP3} where the authors address a problem due to  Muckenhoupt and Wheeden.

\begin{thm}\label{lineargrowth} Let $T$ be a C--Z operator and let $1<p<\infty$. Let $w$ be a weight in $A_1$. Then
\begin{equation*} \label{sharp1}
\|T\|_{L^p(w)}\le c\,pp'\, [w]_{A_1},
\end{equation*}
where $c=c_{n,T}$.
\end{thm}
As an application of this result, the following sharp endpoint estimate can be proven.

\begin{thm}\label{logarithmicgrowth} Let $T$ be a C--Z operator. Let $w$ be a weight in $A_1$. Then
\begin{equation*}\label{improve}
\|T\|_{L^{1,\infty}(w)}\le c\,[w]_{A_1}(1+\log [w]_{A_1}),
\end{equation*}
where $c=c_{n,T}$.
\end{thm}
An important point in the proof of these theorems was the use of a special instance of \eqref{roughC-F} with a gain in the constant $c_{n,w, p}$. To be more precise the following $L^{1}$ estimate was needed:
\begin{equation}
\|Tf\|_{L^1(w)}\le c_{n,T}\, [w]_{A_{q} }\,\|Mf\|_{L^1(w)}   \quad w\in A_{q}
\end{equation}
where there is an improvement in the constant appearing  in \eqref{roughC-F}, because $c_{n,T}$ is a structural constant.

The same kind of qualitative and quantitative Coifman--Fefferman type results are known or can be proved for a variety of singular operators, namely commutators of C-Z operators with BMO functions, vector valued extensions and square functions. In addition, there are also weak type estimates. In that direction, there is also a notion of singularity that can be assigned to these operators. Roughly, these commutators are ``more singular'' than C--Z operators, and these are ``more singular'' than, for example, square functions. This notion of singularity is reflected in the kind of maximal operators involved in the norm inequalities. For the commutator, we have from \cite{CP4} (see also \cite{CP1}) that, for any weight $w\in A_{\infty}$,
\begin{equation*}\label{CFintrobT}
\|[b,T]f\|_{L^p(w)}\le c_{n,w,p}\,\|b\|_{BMO} \|M^2f\|_{L^p(w)},
\end{equation*}
where $M^2=M\circ M$. This result is sharp on the BMO norm of $b$ but not on the weight $w$. More generally, for the iterated commutator $T^k_{b}$, the control operator is $M^{k+1}$. We will show that in this case the constant depends on a $(k+1)$-th power of $[w]_{A_\infty}$. Therefore, as our first main purpose, we will show how to obtain the above mentioned Coifman--Fefferman type inequalities involving Wilson's  $A_{\infty}$ constant $[w]_{A_\infty}$, which is, in light of the chain of inequalities \eqref{mono}, the  good  one.

We will also address the problem of finding sharp operator bounds for a variety of singular operators. The improvement here will be reflected in a \emph{mixed} bound in terms of the $A_{1}$ and $A_{\infty}$  constant of the weight $w\in A_{1}$. This approach is  taken from \cite{HP}, where the authors  prove  a new sharp Reverse H\"older Inequality (RHI) for $A_{\infty}$ weights with the novelty of involving Wilson's constant. It is important to note that this RHI for $A_{q}$ weights, $1\le q\le \infty$, was already known. Moreover, the same proof as in the case of $A_{1}$, which can be found in \cite{LOP1} with minor modifications, also works for $A_{\infty}$ weights but with the ``older'' $\|w\|_{A_{\infty}}$  constant. However, the proof of the same property for  Wilson's constant $[w]_{A_\infty}$ is more difficult and requires a different approach.

Recall that if $w\in A_p$ there are constants $r>1$ and $c\ge1$ such that for any cube $Q$, the RHI holds:
\begin{equation}\label{a2}
\left( \frac{1}{|Q|}\int_Qw^{r}dx\right)^{1/r} \le \frac{c}{|Q|}\int_Qw.
\end{equation}
In the standard proofs both constants $c,r$ depend upon the $A_p$ constant of the weight. A more precise version of \eqref{a2} is the following result that can be found, for instance, in \cite{Pe1}.

\begin{lem}\label{RHAp} Let $w\in A_p$, $1<p<\infty$ and let
$$r_{w}=1+\displaystyle{\frac{1}{2^{2p+n+1}[w]_{A_p}}}.$$
Then for any cube $Q$,
\begin{equation*}\label{a3}
\left( \frac{1}{|Q|}\int_Qw^{r_w}dx\right)^{{1}/{r_w}} \le \frac{2}{|Q|}\int_Qw.
\end{equation*}
\end{lem}

Here we present the new result  of T. Hyt\"onen and the second author (see \cite{HP}). This sharper version of the RHI plays a central role in the proofs of  all the results presented in this article.
\begin{thm}[A new sharp reverse H\"older inequality]\label{thm:SharpRHI}
Define $r_w:=1+\displaystyle{\frac{1}{ \tau_n\,[w]_{A_{\infty}}}}$, where $\tau_n$ is a dimensional constant that we may take to be $\tau_n= 2^{11+n}$.  Note that $r_w'\approx [w]_{A_{\infty}}$.

\vspace{0.2cm}

a)  If $w \in A_{\infty}$, then
$$
\displaystyle\left(\frac{1}{|Q|}\int_Q w^{r_w}\right)^{1/r_w}\leq 2\frac{1}{|Q|}\int_Q w.
$$

b) Furthermore, the result is optimal up to a dimensional factor:  There exists a dimensional constant $c=c_n$ such that, if a weight $w$ satisfies the RHI, i.e., there exists a constant $K$ such that
\begin{equation*}
  \left(\frac{1}{|Q|}\int_Q w^r\right)^{1/r}\leq K\frac{1}{|Q|}\int_Q w,
\end{equation*}
for all cubes $Q$, then $[w]_{A_\infty}\leq  c_n \,K\, r'$.
\end{thm}

Among other important results, in \cite{HP} the authors derive mixed $A_1$--$A_\infty$ type results  of A. Lerner, S. Ombrosi and the second author in \cite{LOP3} of the form:
\begin{equation*}
\|Tf\|_{L^p(w)} \le c\,pp'\, [w]_{A_1}^{1/p}[w]_{A_{\infty}}^{1/p'}, \qquad  w\in A_{1}, \quad 1<p<\infty,
\end{equation*}
which is an improvement of Theorem \ref{lineargrowth}. They also derive weak type estimates like
\begin{equation*}
\|Tf\|_{L^{1,\infty}(w)} \le c[w]_{A_1} \log(e+[w]_{A_{\infty}})\|f\|_{L^1(w)}.
\end{equation*}

As another main purpose in the present work, we will show how to extend this results on mixed $A_1$--$A_\infty$ bounds to the case of the commutator and its higher order analogue. The analogues of Theorem \ref{lineargrowth} and  Theorem \ref{logarithmicgrowth}  for the commutators we consider were already proved by the first author in \cite{O}:

\begin{thm}[Quadratic $A_{1}$ bound for commutators]\label{quadraticA1_bT}
Let $T$ be a Calder\'{o}n--Zygmund operator and let $b$ be in $BMO$.
Also let $1<p,r<\infty$. Then there exists a constant $c=c_{n,T}$
such that for any weight $w$, the following inequality
holds
\begin{equation}\label{acotadospesos}
\|[b,T]f\|_{L^{p}(w)}\le c\,  \|b\|_{BMO}\,
{(pp')}^2\,(r')^{1+\frac{1}{p'}}\,
   \|f\|_{L^{p}(M_{r}w)}.
\end{equation}
In particular if $w \in A_1$, we have
\begin{equation}
\|[b,T]f\|_{L^p(w)}\le c\,\|b\|_{BMO} (pp')^2[w]_{A_1}^2\,
   \|f\|_{L^{p}(M_{r}w)}.
\end{equation}
Furthermore this result is sharp in both $p$ and the exponent of  $[w]_{A_1}$.
\end{thm}

\begin{thm}\label{debil}
Let $T$ and $b$ as above. Then there exists a constant $c=c_{n,T}$
such that for any weight

\begin{equation}\label{comdebil1}
w(\{x\in \mathbb R^n:|[b,T]f(x)|>\lambda\})\leq c\, (pp')^{2p}(r')^{2p-1}\int_{\mathbb R^n} \Phi
\left(\|b\|_{BMO}{\frac{|f|}{\lambda}} \right)\,M_rwdx,
\end{equation}
where $\Phi(t)=t(1+\log^{+}t)$.

As a consequence, If $w \in A_1$
\begin{equation}
w(\{x\in \mathbb R^n:|[b,T]f(x)|>\lambda\})\leq c\,
\Phi([w]_{A_1})^2 \int_{\mathbb R^n} \Phi
\left({\frac{|f(x)|}{\lambda}} \right)w(x)\,dx
\end{equation}
where $\Phi(t)=t(1+\log^{+}t)$.
\end{thm}

In this article we present an improvement of these theorems in terms of mixed $A_{1}$--$A_{\infty}$ norms for the commutator and for its iterations, proved by the first author in her dissertation (2011). For any $k\in\mathbb{N}$, the $k$-th iterated commutator $T^{k}_{b}$ of a BMO function $b$ and a C-Z operator $T$ is defined by
$$ T^k_{b}:=[b,T_{b}^{k-1}].$$

This paper is organized as follows. In Section \ref{main-results} we formulate the precise statements of the results announced in the introduction. In Section \ref{prelim} we present some background and auxiliary results for the proofs and finally, in Section \ref{proofs}, we describe the main features of the proofs.

\section{Main Results}\label{main-results}

In this section we will present the precise statements of the results discussed in the introduction. First, we present  several Coifman--Fefferman inequalities for a variety of singular operators. Our purpose is to emphasize that there is a notion of order of singularity that allows us to distinguish them. This higher or lower singularity can be seen in the power of the maximal functions involved and also in the dependence of the  constant of the weight. Next, we present our results on strong and weak norm inequalities for commutator and for its  iterations. In this latter case, we study in detail how the iterations affects on each part of the weight, namely the $A_{1}$ and the $A_{\infty}$ fractions of the constant.

\subsection{Coifman-Fefferman inequalities}\label{CF-mainresults}

\begin{thm}[C-Z operators]\label{CZ}
Let $T$ be a C--Z operator and let $w\in A_{\infty}$. Then there is constant $c=c_{n,T}$ such that, for any $0<p<\infty$,
\begin{equation*}
\|Tf\|_{L^p(w)}\le c\,\max\{1,p\}[w]_{A_{\infty}}\|Mf\|_{L^p(w)}
\end{equation*}
whenever $f$ is a function satisfying the condition $|\{ x : |Tf(x)|> t\}| <\infty$ for all $t>0$.
\end{thm}

We need some additional notation for the following theorem. Let $1<q<\infty$ and define, for  any sequence of functions $f=\{f_j\}_{j\in\mathbb{N}}$, $|f|_q:= \left(\sum_j |f_j|^q\right)^\frac1{q}$. Also define, for a C-Z oprator $T$, the vector valued extension $\overline{T}f=\{Tf_j\}_{j\in\mathbb{N}}$.

\begin{thm}[C-Z operators - Vector valued extensions]\label{CZVE}
Let $T$ be a C--Z operator and let $w\in A_{\infty}$. Then there is constant $c=c_{n,T}$ such that, for any $0<p<\infty$, $1<q<\infty$ and for any sequence of compactly supported functions $\{f_{j}\}_{j\in\mathbb{N}}$,
\begin{equation*}
\left\||\overline{T}f|_{q}\right\|_{L^p(w)}\le c\max\{1,p\}[w]_{A_{\infty}}\left\|M\left(|f|_{q}\right)\right\|_{L^p(w)}.
\end{equation*}
\end{thm}

We also have the following result on multilinear C--Z operators. For the precise definitions and properties, see \cite{LOPTTG}. Let $T$ be an $m$-linear C--Z operator acting on a vector $\vec{f}\,$ of $m$ functions $\vec{f}\,=(f_{1},\dots,f_{m})$. Define also the maximal function $\mathcal{M}$ by
	$$\mathcal{M}(\vec{f}\,)(x)=\sup_{Q\ni x} \prod_{i=1}^m\frac{1}{|Q|}\int_{Q}|f_{i}(y_{i})|\ dy_{i}.$$
The following theorem is a refinement of \cite[Corollary 3.8]{LOPTTG}.

\begin{thm}[Multilinear C--Z operators]\label{MCZ}
Let $T$ be an $m$-linear C--Z operator, let $w\in A_{\infty}$ and let $p>0$. Then there exists a constant $c=c_{n,m,T}$ such that
\begin{equation*}
\|T(\vec{f}\,)\|_{L^p(w)}\le c\,\max\{1,p\}[w]_{A_{\infty}}\|\mathcal{M}(\vec{f}\,)\|_{L^p(w)},
\end{equation*}
whenever $\vec{f}\,$ is a vector of compactly supported functions.
\end{thm}

\begin{thm}[Commutators]\label{Comm}
Let $T$ be a C--Z operator, $b\in BMO$ and let $w\in A_{\infty}$. Then there is a constant $c=c_{n,T}$ such that for $0<p<\infty$,
\begin{equation*}
 \|[b,T]\|_{L^p(w)}\le c\,\|b\|_{BMO}\max\{1,p^2\}  [w]_{A_\infty}^2 \|M^2f\|_{L^p(w)}
\end{equation*}
\end{thm}

More generally, we have  the following result for the iterated commutator.
\begin{thm}[$k$-th Iterated Commutator]\label{itComm}
Let $T$ be a C--Z operator, $b\in BMO$ and let $w\in A_{\infty}$. For the $k$-th $(k\ge2)$ commutator $T^k_{b}$, there is a constant $c=c_{n,T}$ such that, for $0<p<\infty$,
\begin{equation*}
 \|T^k_{b}f\|_{L^p(w)}\le c\, 2^k\max\{1,p^{k+1}\}  [w]_{A_\infty}^{k+1} \|b\|^k_{BMO}\|M^{k+1}f\|_{L^p(w)}
\end{equation*}
\end{thm}

\subsection{Mixed \texorpdfstring{$A_{1}$--$A_{\infty}$}{a1-ainf} strong and weak norm inequalities for commutators}\label{mixed-mainresults}

\begin{thm}
Let $T$ be a  C--Z operator, $b\in BMO$ and let $1<p<\infty$.  Then there is a constant $c=c_{n,T}$   such that, for any $w\in A_1$,
\begin{equation*}
\|{[b,T]}\|_{L^p(w)}\le c\,\|{b}\|_{BMO} (pp')^2[w]_{A_1}^{1/p}[w]_{A_{\infty}}^{1+1/p'}.
\end{equation*}
\end{thm}

\begin{thm}
Let $T$ and $b$ as above. Then there exists a constant $c=c_{n,T}$ such that for any weight $w \in A_1$
\begin{equation*}
w(\{x\in \mathbb R^n:|[b,T]f|>\lambda\})\leq c\,
\beta \int_{\mathbb R^n} \Phi
\left(\|{b}\|_{BMO}{\frac{|f|}{\lambda}} \right)\,w(x)dx,
\end{equation*}
where
$\beta= [w]_{A_1}[w]_{A_{\infty}}(1+\log^{+}[w]_{A_{\infty}})^2$ and $\Phi(t)=t(1+\log^{+}t)$.
\end{thm}

More generally, we have the following generalization for the $k$-th iterated commutator.

\begin{thm}
Let $T$ be a C--Z operator, $b\in BMO$ and  $1<p,r<\infty$. Consider the higher order commutators $T^k_{b}$,  $k=1,2,\cdots,$. Then there exists a constant $c=c_{n,T}$ such that for any weight $w$ the following inequality holds
\begin{equation}\label{kacotadospesos}
\|T_{b}^{k}f\|_{L^{p}(w)}\le c\,  \|b\|_{BMO}^{k}\,
{(pp')}^{k+1}\,(r')^{k+{1/p'}}\,
   \|f\|_{L^{p}(M_{r}w)}.
\end{equation}

In particular if $w\in A_1$, we have that
\begin{equation*}
\|{T_{b}^{k}}\|_{L^p(w)}\le c\,\|{b}\|_{BMO}^k (pp')^{k+1}[w]_{A_1}^{1/p}[w]_{A_{\infty}}^{k+1/p'}.
\end{equation*}
\end{thm}

\begin{thm}
Let $T$ and $b$ as above, and let $1<p,r<\infty$. Then there exists a constant $c=c_{n,T}$
such that for any $w$
\begin{equation}\label{kcomdebil1}
w(\{x\in \mathbb R^n:|T_{b}^{k}f|>\lambda\})\leq c\,
(pp')^{(k+1)p}(r')^{(k+1)p-1} \int_{\mathbb R^n} \Phi
\left(\|b\|_{BMO}{\frac{|f|}{\lambda}}\right) M_rw\,dx,
\end{equation}
where $\Phi(t)=t(1+\log^{+}t)^k$.

If $w\in A_1$ we obtain that
\begin{equation*}
w(\{x\in \mathbb R^n:|T_b^k f(x)|>\lambda\})\leq c_n\,
\beta \,\int_{\mathbb R^n} \Phi
\left(\|{b}\|_{BMO}{\frac{|f(x)|}{\lambda}} \right)\,w(x)dx,
\end{equation*}
where $\beta= [w]_{A_1}[w]^k_{A_{\infty}}(1+\log^{+}[w]_{A_{\infty}})^{k+1}$ and $\Phi(t)=t(1+\log^{+}t)^k$.
\end{thm}

\section{Background and Preliminaries} \label{prelim}

\subsection{Rearrangement type estimates}

In this section we present an important tool based on rearrangements of functions. The main lemma is the following about the control of $L^p$ norms by a sharp maximal function. We start by recalling some standard definitions. Given a locally integrable function $f$ on $\mathbb{R}^n$, the \emph{Hardy--Littlewood maximal operator} $M$ is defined by
\begin{equation*}
Mf(x)=\sup_{Q\ni x}{\frac1{|Q|}}\int_{Q} f(y)\\,dy,
\end{equation*}
where the supremum is taken over all cubes $Q$ containing the point $x$. We will also use the following operator:
$$M_\delta f(x)=(M(|f|^\delta)(x))^{1/\delta}, \qquad 0<\delta < 1.$$
Recall that the Fefferman-Stein sharp maximal function is defined as

$$M^\#(f)(x)=\sup_{Q\ni x}\frac{1}{|Q|}\int_Q |f(y)-f_Q|\, dy,$$
where

$$ f_Q=\frac{1}{|Q|}\int_Q f(y)\, dy.$$

If we  only consider dyadic cubes, we obtain the \emph{dyadic} sharp maximal function, denoted by $M^{\#,d}$.

Define also, for $0<\delta<1$,

$$M_\delta ^\# f(x)=M^\# (|f|^\delta)(x)^{1/\delta}.$$

and

$$M_\delta ^{\#,d} f(x)=M^{\#,d} (|f|^\delta)(x)^{1/\delta}.$$

\begin{lem}\label{Bagby-Kurtz}
Let $0<p<\infty$, $0<\delta<1$ and let $w\in A_\infty$. Then
\begin{equation}\label{Fefferma-Stein1}
\|f\|_{L^p(w)}\le c\max\{1,p\}[w]_{A_\infty}\,\|M_{\delta}^{\#,d}(f)\|_{L^p(w)}
\end{equation}

for any function $f$ such that $|\{ x : |f(x)|> t\}| <\infty$ for all $t>0$.
\end{lem}

The proof of Lemma \ref{Bagby-Kurtz} will follow from an analogous inequality for the non-increasing rearrangements of $f$ and $M^{\#,d}_{\delta}$ with respect to the weight $w$. Recall that the non-increasing rearrangement $f_{w}^*$ of a measurable function $f$ with respect to a weight $w$ is defined by
$$ f_w^*(t) = \inf\big\{\lambda>0:\, w_f(\lambda)< t\big\},\qquad t>0,$$
where
$$w_f(\lambda)= w(\{x\in \mathbb{R}^n:\,|f(x)|>\lambda\}),\,\,\lambda>0,$$
is the distribution function of $f$ associated to $w$. An important fact is that
$$\int_{\mathbb{R}^n}|f|^p\,wdx = \int_0^{\infty} f_w^*(t)^p\,dt.$$

The key rearrangement lemma is the following.
\begin{lem}\label{keyrearreng}
Let and $w\in A_{\infty}$, $0<\delta<1$ and $0<\gamma<1$. There is a
constant $c=c_{n,\gamma,\delta}$, such that for any measurable function:
\begin{equation} \label{keyrearreng1}
f^*_{w}(t) \le c[w]_{A_\infty} \,\big(M_{\delta}^{\#,d}f\big)^*_{w}(\gamma\,t)+
f^*_{w}(2t) \qquad t>0
\end{equation}
\end{lem}
These sort of estimates goes back to the work of R. Bagby and D. Kurtz in the middle of the  80's (see \cite{BK1} and \cite{BK2}). The proof of Lemma \ref{keyrearreng} can be found in \cite{Pe1} but with $A_{q}$ instead of $A_{\infty}$  weights.  This improvement to $A_{\infty}$ weights and, moreover, with the smaller $[w]_{\infty}$ constant, is a consequence of the sharp RHI for $A_{\infty}$ weights from \cite{HP}. We will sketch the proof in the next section. With this lemma, we can prove Lemma \ref{Bagby-Kurtz}. If we iterate \eqref{keyrearreng1} we have:
\begin{eqnarray*}
f^*_{w}(t) & \leq  & c\,[w]_{A_\infty}\,\sum^{\infty}_{k=0}(M_{\delta}^{\#}f)^*_{w}(2^{k}\gamma t)+ f^*_{w}(+\infty)\\
& \le & \frac{[w]_{A_\infty}\,}{\log 2} \,\int_{t\gamma/2}^{\infty}(M_{\delta}^{\#}f)^*_{w}(s)\,\frac{ds}{s}+ f^*_{w}(+\infty).
\end{eqnarray*}
Hence if we assume that
\begin{equation*}\label{assump.onf}
  f^*_{w}(+\infty)=0,
\end{equation*}
the inequality we obtain is
\begin{equation*}\label{PointwRearr.}
f^*_{w}(t) \leq c\,[w]_{A_\infty} \,\int_{t\gamma/2}^{\infty}
(M_{\delta}^{\#}f)^*_{w}(s)\,\frac{ds}{s}.
\end{equation*}

We continue using the Hardy operator. Recall that if $f:(0,\infty)\rightarrow [0,\infty)$
$$
Af(x)=\frac1x \int_0^x f(t)dt, \qquad x>0
$$
is called the Hardy operator. The dual operator is given by
$$
Sf(x)= \int_x^{\infty} f(s) \frac{ds}{s}.
$$

Hence the above estimate can be expressed as:
$$ f^*_{w}(t) \leq c\, [w]_{A_\infty} \,S((M_{\delta}^{\#}f)^*_{w})(t\gamma/2)$$

Finally since it is well known that these operators are bounded on
$L^p(0,\infty)$ and furthermore, it is known that if $p\geq 1$, then $\|S\|_{L^p(0,\infty)}= p$,
we have that
\begin{eqnarray*}
\|f\|_{L^p(w)}&=& \|f_w^*\|_{L^p(0,\infty)}\\
 & \le & c\,[w]_{A_\infty}\, \|S((M_{\delta}^{\#}f)_w^*)\|_{L^p(0,\infty)}\\
& \le & c\,p\,[w]_{A_\infty}\, \|(M_{\delta}^{\#}f)_w^*\|_{L^p(0,\infty)}\\
& = & c\,p\,[w]_{A_\infty}\,\|M_{\delta}^{\#}f\|_{L^p(w)}
\end{eqnarray*}

This concludes the strong estimate in the case $p\geq 1$. For $0<p<1$, use the triangle inequality and $L^1$ boundedness of the operator $S$ (we get no $p$ here, so we just get $\max\{1,p\}$ for general $p$). This concludes the proof of Lemma \ref{Bagby-Kurtz}.

\subsection{Pointwise inequalities}

Once we have the result of the previous subsection, i.e. the Fefferman--Stein inequality involving $M^{\#,d}_{\delta}$ with a sharp dependence on the weight, we can then deduce sharp Coifman--Fefferman weighted inequalities for any pair of operators $T$ and $S$ satisfying a pointwise inequality like
$$ M^{\#}_{\delta}(Tf)(x)\le c_{\delta} Sf(x)\qquad \text{a.e.} \ x\in\mathbb{R}^n.$$

We will elaborate on this in the next section. Here we want to collect those inequalities that involve the operators under study (we refer the interested reader to  \cite[Chapter 9]{CMP3} for a deeper treatment of this subject). The inequalities that are going to be used are:

\begin{itemize}
\item For $0<\delta<1$ and $T$ any C--Z operator (see \cite{AP}),
\begin{equation}\label{maximal(T)vsMaximal}
M^\#_{\delta}(Tf)(x)\le c_{\delta}\,Mf(x).
\end{equation}
\item We have the following vector valued extension from \cite{PT}: Let $1<q<\infty$ and $0<\delta<1$. Let $T$ be a C-Z operator and consider the vector valued extension $\overline{T}$. There exists a constant $c_{\delta}$ such that
	$$ M^\#_{\delta}\left(|\overline{T}f|_q\right)(x)\le c_{\delta}\,M(|f|_q)(x).$$

\item We also have a pointwise inequality for multilinear C--Z operators. Let $T$ be an $m$-linear C--Z operator. Then for all $\vec{f}\,$ in any product of $L^{q_{j}}(\mathbb{R}^n)$ spaces, with $1\le q_{j}<\infty$,
	$$ M^\#_{\delta}(T(\vec{f}\,))(x)\le c_{\delta}\,\mathcal{M}(\vec{f}\,)(x)\qquad  \text{ for } 0<\delta<1/m.$$

\item The following result is from \cite{O}. If $0<\delta<\varepsilon<1$, there is a constant $c=c_{\varepsilon,\delta}$ such that
	$$M^{\#,d}_{\delta}(M^d_{\varepsilon}(f))(x)\le c \,M^{\#,d}_{\varepsilon}f(x).$$
	What we really need,  for the proof of the Coifman--Fefferman inequality in the case of the commutator, is a consequence of this inequality, which is a subtle improvement of Lemma \ref{Bagby-Kurtz}: for $0<p<\infty$ and $0<\delta<1$,
	\begin{equation}\label{maximal-Bahby-Kurtz}
	\|M^d_{\delta}f\|_{L^p(w)}\le c\,\max\{1,p\}[w]_{A_\infty}\,\|M_{\delta}^{\#,d}(f)\|_{L^p(w)}.
	\end{equation}
\item Another related and very important pointwise inequality for iterated commutators is the following result from \cite{CP4}. For each $b\in BMO$, $0<\delta<\varepsilon<1$, there exists $c=c_{\delta,\varepsilon}$ such that, for all smooth functions $f$, we have that
\begin{equation}\label{pointComm}
M^\#_{\delta}(T^k_{b}f)(x)\le c\,\|b\|_{BMO}\sum_{j=0}^{k-1}M^d_{\varepsilon}(T^j_{b}f)(x) + \|b\|_{BMO}^kM^{k+1}f(x).
\end{equation}

\end{itemize}

\subsection{Building \texorpdfstring{$A_{1}$}{a1} weights from duality.}
The following lemma gives a way to produce $A_{1}$ weights with special control on the constant. It is based on the  so called Rubio de Francia iteration scheme or algorithm.

\begin{lem}\cite{LOP3}\label{rubiodefrancia}
Let $1<s<\infty$, and let $v$ be a weight. Then, there exists a
nonnegative sublinear operator $R$ bounded in $L^s(v)$ satisfying the following
properties:
\begin{enumerate}
  \item[(i)] $h\leq R(h)$;
  \item[(ii)] $\|{Rh}\|_{L^{s}(v)}\leq 2\|{h}\|_{L^{s}(v)}$;
  \item[(iii)] ${Rh\,v^{1/s}\in A_1}$ with
  \begin{equation*}
 [Rh\,v^{1/s}]_{A_1}\leq cs'.
  \end{equation*}
\end{enumerate}
\end{lem}

\section{About the proofs}\label{proofs}

\subsection{Coifman--Fefferman inequalities}

We begin this section with a brief description of how to obtain Coifman-Fefferman inequalities for all the operators presented in the previous section. These inequalities will follow from a general scheme involving pointwise control of the singular operator by an appropriate maximal function and the key inequality from Lemma \ref{Bagby-Kurtz}. Suppose that we have, for a pair of operators $T$ and $S$, a pointwise inequality like
$$ M^{\#}_{\delta}(Tf)(x)\le c_{\delta}\, Sf(x)\qquad \text{a.e.} \ x\in\mathbb{R}^n.$$
Then, for any $f$ such that  $|\{ x : |Tf(x)|> t\}| <\infty$ for all $t>0$ and any $w\in A_{\infty}$, we have that
\begin{eqnarray*}
\|T(f)\|_{L^p(w)} & \le & c\,\max\{1,p\}[w]_{A_\infty}\,\|M_{\delta}^{\#,d}(T(f))\|_{L^p(w)}\\
& \le & c_{\delta}\,\max\{1,p\}[w]_{A_\infty}\,\|S(f)\|_{L^p(w)}
\end{eqnarray*}
This argument completes the proof of Theorem \ref{CZ}, Theorem \ref{CZVE} and Theorem  \ref{MCZ}.

For the case of the commutator and its iterations, we need an extra step for each iteration. More precisely, for the commutator we have the following estimate. Let $w$ be  any $A_{\infty}$ weight and let $p>0$.  For $0<\varepsilon<\delta<1$, we have that
\begin{eqnarray*}
\|[b,T]f\|_{L^p(w)} &\le &  c\,\max\{1,p\}[w]_{A_\infty}\,\|M_{\varepsilon}^{\#,d}([b,T]f)\|_{L^p(w)}\\
& \le & c\,\max\{1,p\}[w]_{A_\infty}\|b\|_{BMO}\left(\|M_{\delta}^{d}(Tf)\|_{L^p(w)}+ \|M^2f\|_{L^p(w)}\right)\\
\end{eqnarray*}
We start with \eqref{Bagby-Kurtz} applied to $[b,T]f$, and then we use \eqref{pointComm}. Now we can combine \eqref{maximal-Bahby-Kurtz} with \eqref{maximal(T)vsMaximal} to control the first term and therefore obtain that
\begin{equation*}
\|[b,T]f\|_{L^p(w)} \le    c\, \|b\|_{BMO}\max\{1,p^2\}[w]^2_{A_\infty}\|b\|_{BMO}\|M^2f\|_{L^p(w)}.
\end{equation*}
This final estimate follows from the fact that $[w]_{A_{\infty}}\ge 1$ and by dominating $M$ by $M^2$. For the iterated commutator, we proceed in the same way as before. The relevant inequality, which can be easily proved by induction, is contained in the following lemma.
\begin{lem}
Let $0<\delta<1$. Then there exists $\eta$ such that $0<\delta<\eta<1$ and a universal constant  $c>0$ such that
$$ \|M_{\varepsilon}^{\#,d}(T_{b}^kf)\|_{L^p(w)}\le 2^k\, c\,\max\{1,p\}\|b\|_{BMO}\|M^{k+1}_{\eta}f\|_{L^p(w)}$$

\end{lem}

Following the preceding scheme and applying this lemma, we obtain the result stated in Theorem \ref{itComm}. This will conclude the proofs of all the announced Coifman--Fefferman inequalities, if we prove the essential inequality \eqref{keyrearreng1}. In fact, it is in this inequality where the improvement on the dependence on the weight appears. We will present here a sketch of the proof from \cite{Pe1}. The idea is to show that the key ingredient for the improvement is an exponential decay lemma and the new Reverse H\"older Property for $A_{\infty}$ weights.

\subsubsection{Proof of the key rearrangement estimate}

Fix $t>0$ and let $A=[w]_{A_\infty}$. By definition of rearrangement, \eqref{keyrearreng1} will follow if we prove that for each $t>0$,

\begin{equation*}
w\Big\{ x \in \mathbb{R}^n  : f(x)>
cA\,\big(M_{\delta}^{\#}f\big)^*_{w}(\gamma\,t)+ f^*_{w}(2t)  \Big\} \le t
\end{equation*}

We split the left hand side $L$ as follows:

\begin{eqnarray*}
L&\le&  w\{ x \in \mathbb{R}^n:\,  cA\,M_{\delta}^{\#}f(x)> cA\,\big(M_{\delta}^{\#}f\big)^*_{w}( \gamma\,t)\\
&+& w\{ x \in \mathbb{R}^n : f(x)> cA\,M_{\delta}^{\#}f(x)+ f^*_{w}(2t) \}\\
&=& I+II.
\end{eqnarray*}
Observe that
$$I= w\Big\{ x \in \mathbb{R}^n : \,M_{\delta}^{\#}f(x)> \big(M_{\delta}^{\#}f\big)^*_{w}(
\gamma\,t)   \Big\}\leq \gamma\,t$$
by definition of rearrangement, and hence the heart of the matter is to prove that
$$II=w\{ x \in \mathbb{R}^n : f(x)> cA\,M_{\delta}^{\#}f(x)+ f^*_{w}(2t)
\}\le (1-\gamma)\,t$$

The set $\{ x \in \mathbb{R}^n : f(x)> cA\,M_{\delta}^{\#}f(x)+ f^*_{w}(2t)\}$ is contained in the set $E=\{x \in \mathbb{R}^n: f(x)> f^*_{w}(2t)\}$ which has $w$--measure at most $2t$ by definition of rearrangement. Now, by the regularity of the measure we can find an open set $\Omega$ containing $E$ such that $w(\Omega)<3t$. The remainder of the proof from \cite{Pe1} consists in proving that there exists a constant $c$ depending only on $\delta$ and the dimension such that
$$II= w\{ x \in \Omega : f(x)> cA\,M_{\delta}^{\#}f(x)+ f^*_{w}(2t)\} \le (1-\gamma)\,t.$$
By applying an appropriate Calder\'on-Zygmund decomposition of the set $\Omega$ into a union of dyadic cubes $\{Q_{j}\}$, it can be proved that, for $c>1+2^{1/\delta}$, we have the following inclusion for any $j$. If we denote
\begin{equation}\label{intialpoint}
E_j= \{x\in Q_j: |f(x)-m_f(Q_j)|> cA\,M_{\delta}^{\#}f(x) \Big\},
\end{equation}
then
$$
\{x\in Q_j: f(x)> cA\,M_{\delta}^{\#}f(x)+ f^*_{w}(2t)   \} \subset E_{j},
$$
where $m_f(Q)$ is the median value of $f$ over $Q$ (see \cite{L04}). The key ingredient now is an ``exponential decay'' Lemma from \cite{Pe1} adapted to $E_{j}$.
\begin{lem}\label{expdecay}
Let $f \in L^{\delta}_{loc}$. For $t>0$ we define
$$
\varphi(t) =  \sup_{Q \in \mathcal{D} }  \frac{1}{|Q|} \Big|\{x\in
Q: |f(x)-m_f(Q)|
>t\,M_{\delta}^{\#}f(x)\} \Big|
$$
There are dimensional constants $c_1,c_2$ such that $\varphi(t) \leq {\displaystyle\frac{c_1}{e^{c_2t}}}.$
\end{lem}
Applying this lemma to \eqref{intialpoint}, we obtain that
$$ \frac{|E_j|}{|Q_j|}\leq c_1e^{-c_2\, cA}.$$
Now we apply H\"older's inequality to obtain that
\begin{eqnarray}\label{HolderTrick}
\nonumber w(E_{j})& = & \frac{1}{|Q_{j}|}\int_{Q_{j}}\chi_{E_{j}}(x)w(x)\ dx\ |Q_{j}|\\
& \le & \left(\frac{|E_{j}|}{|Q_{j}|}\right)^{1/r'}\left(\frac{1}{|Q_{j}|}\int_{Q_{j}}w(x)^r\ dx\right)^{1/r}|Q_{j}|
\end{eqnarray}
We conclude by choosing $r=r_w$, the sharp reverse exponent for $w$, (see Theorem \ref{thm:SharpRHI}), $r_{w}=1+\displaystyle{\frac{1}{2^{11+d}[w]_{{A_{\infty}}}}}$ . We obtain
$$
\left( \frac{1}{|Q_{j}|}\int_{Q_{j}}w^{r_w}dx\right)^{1/r_w} \le
\frac{2}{|Q_{j}|}\int_{Q_{j}}w
$$
and therefore

$$ w(E_j) \le 2\,\left(\frac{|E_j|}{|Q_j|}\right)^{1/r'} w(Q_j) \leq c_{1}e^{-c\,c_2} w(Q_j)
$$
since $r'\approx [w]_{A_\infty}=A$. Here $c$ and $c_2$ depend on the dimension and we still can choose $c$ as big as needed. If we finally we choose $c$ such that $c_{1}{e^{-c_2\,c}}<\frac{1-\gamma}{3}$, then we have that
\[
II \le \sum_{j} w(E_j)\le  \frac{1-\gamma}{3} \sum_{j} w(Q_j) =
\frac{1-\gamma}{3}\,w(\Omega) < (1-\gamma)\,t,
\]
since $w(\Omega)< 3t$.

As a final remark on Coifman-Fefferman type inequalities, we want to mention that  for the case of C-Z operators, Theorem \ref{CZ} can be deduced from a combination of standard good-$\lambda$ techniques together with the following local exponential estimate from \cite{OPR}, which is an improvement of the classical result of Buckley \cite{B} (see also \cite{K}).

\begin{thm}\label{thm:subexpT}
Let $T$ be a Calderon-Zygmund operator. Let $Q$ be a cube and let $f\in L_{c}^\infty(\mathbb{R}^n)$  such that $\emph{supp}(f)=Q$. Then there exist $c>0$ and $k>0$ such that
$$|\{x\in Q: |Tf(x)|>tMf(x) \}|\le k e^{-ct}|Q|,\qquad t>0.$$
\end{thm}
This theorem is essentially all we need for the proof of Theorem \ref{CZ} since, by standard truncation arguments, we can restrict ourselves to look at only the local part. Now, if we put $E_j:=\{x\in Q_{j}:T(f)>\lambda, Mf<\gamma \lambda\}$, we have that,  for some constant $c_{1}$,
\begin{equation*}
\frac{|E_{j}|}{|Q_{j}|}\le c_1e^{-c_2/\gamma}, \qquad \lambda>0,\gamma_{0}>\gamma>0.
\end{equation*}
Here we are in the exact same setting as in \eqref{HolderTrick}, and therefore we can obtain that
$$ w(E_j) \le 2\,\left(\frac{|E_j|}{|Q_j|}\right)^{1/r'} w(Q_j) \leq c_{1}e^{-c_2/\gamma r' } w(Q_j),$$
where $r'\approx [w]_{A_\infty}$. Now we sum over $j$ to obtain for any constant $B>0$ the following good-$\lambda$ inequality:
\begin{equation*}
w\left\{x\in\mathbb{R}^n:T(f)>3\lambda, Mf<\frac{B}{[w]_{\infty}}\lambda\right\}\le c_{1}e^{-{c_{2}/B}}w\left\{x\in\mathbb{R}^n:T(f)>\lambda\right\}
\end{equation*}
From this inequality, we can apply standard good-$\lambda$ arguments to obtain a sharp Coifman--Fefferman inequality with $[w]_{\infty}$ bound.

We also remark that in \cite{OPR} a similar local \emph{sub}exponential decay is proved for the commutator $[b,T]$ above. It is natural to ask if it is possible to prove a Coifman--Fefferman inequality for $[b,T]$ via truncation, exponential decay and sharp RHI. Unfortunately, we do not know how to control the part of the truncated operator ``away'' from some fixed cube.

\subsection{Mixed \texorpdfstring{$A_{1}$--$A_{\infty}$}{a1-ainf-comm} strong and weak norm inequalities for commutators}

In this section we present the proofs of the results announced in Section \ref{mixed-mainresults}. Those results actually follow from the two weight bounds for commutators \eqref{acotadospesos} proved in \cite{O} for $k=1$, and similarly \eqref{kacotadospesos} (see\cite{O2}).

Let us briefly describe how to derive this inequality in the case $k= 1$. The aim is to prove that
\begin{equation}
\|{[b,T]f}\|_{L^{p}(w)}\le c\,  \|{b}\|_{BMO}\,
{(pp')}^{2}\,(r')^{1+\frac{1}{p'}}\,
   \|{f}\|_{L^{p}(M_{r}w)}.
\end{equation}
By duality, it is equivalent to prove that
\begin{equation*}\label{kacotadospesosdual}
\left\|{\displaystyle{\frac{([b,T])^*f}{M_{r}w}}}\right\|_{L^{p'}(M_{r}w)}\le
c\, (pp')^{2}\,
(r')^{1+\frac{1}{p'}}\,\left\|{\displaystyle{\frac{f}{w}}}\right\|_{L^{p'}(w)},
\end{equation*}
for the adjoint operator $[b,T]^*$. To prove this last inequality, we write the norm as
\begin{equation*}
\left\|{\displaystyle{\frac{([b,T])^*f}{M_{r}w}}}\right\|_{L^{p'}(M_{r}w)}=\sup_{\displaystyle{\|{h}\|_{L^{p}(M_{r}w)}=1}}
\left|\int_{\mathbb{R}^n}([b,T])^*f(x)h(x)\, dx \right|.
\end{equation*}
Now we use Rubio de Francia's algorithm and Lemma \ref{Bagby-Kurtz} to obtain that
\begin{eqnarray*}
\left|\int_{\mathbb{R}^n}([b,T])^*f(x)h(x)\, dx \right|\leq
c_{\delta}\,[Rh]_{A_3}\int_{\mathbb{R}^n}M_{\delta}^{\#}(([b,T])^*f)(x)\,Rh(x)\,dx\\
\end{eqnarray*}
The key here is to use that $[Rh]_{A_3}\leq c_n\,p'$. Also note that for our purposes here it is sufficient to use  \eqref{Fefferma-Stein1} with the $A_{3}$ constant of the weight. To handle the integral we use the pointwise inequality \eqref{pointComm}, and therefore we obtain that
\begin{equation*}
\left|\int_{\mathbb{R}^n}([b,T])^*f(x)h(x)\, dx \right|\leq
c_{\delta,n,\varepsilon}\,p'\|b\|_{BMO}\int_{\mathbb{R}^n}M_{\varepsilon}^d(T^*f)(x)+M^2f\,Rh(x)\,dx.\\
\end{equation*}
For the second term we only need to use H\"older and the properties of the weight $Rh$. For the first term, we use \eqref{maximal-Bahby-Kurtz} and then \eqref{maximal(T)vsMaximal}. We obtain that
\begin{equation*}
\left|\int_{\mathbb{R}^n}([b,T])^*f(x)h(x)\, dx \right|\leq
c_{\delta,n,\varepsilon}\,(p')^2\|b\|_{BMO}\|M^2f\|_{L^{p'}((M_{r}w)^{1-p'})}.
\end{equation*}
Finally, we use  (see \cite{O})  that
$$\|M^2f\|_{L^{p'}((M_{r}w)^{1-p'})}\le c p^2(r')^{1+1/p'}\|f\|_{L^p(w^{1-p'})},$$
which implies the desired result.

The  case $k>1$ is technically more difficult, but the ideas are similar. The key estimate is the pointwise inequality \eqref{pointComm} and  an induction argument (see \cite{O2} for details).

Once we have proved \eqref{kacotadospesos}, we can derive the announced mixed bound as follows.
First, note that it is easy to verify that for \emph{any} $r>1$,
\begin{equation*}\label{MaxrVSw}
M_{r}w(x)\le 2[w]_{A_{1}}w(x),
\end{equation*}
which implies in a standard way a weighted estimate with sharp dependence on $[w]_{A_{1}}$. But we are able to obtain an  improvement by choosing the best RHI at our disposal. Choose then $r=r_{w}$ as the sharp exponent in the RHI for $A_{\infty}$ weights, i.e. $r=r_{w}=1+\displaystyle{\frac{1}{c_{n}[w]_{A_{\infty}}}}$. In this case we get a fraction of the weight from the pointwise inequality \eqref{MaxrVSw} applied to the integrand and, in addition, now $r_{w}'$ is comparable with $[w]_{A_{\infty}}$, and therefore from \eqref{kacotadospesos} we conclude that
\begin{equation*}
\|{T_{b}^{k}f}\|_{L^{p}(w)}\le c\,  \|{b}\|_{BMO}^{k}\, {(pp')}^{k+1}\,[w]_{A_{1}}^{\frac{1}{p}} [w]_{A_{\infty}}^{k+\frac{1}{p'}}\,\|{f}\|_{L^{p}(w)}.
\end{equation*}

\bigskip

Now we present the weak-type inequalities. In the same way as before, the mixed bounds will follow from  a two weight result \eqref{comdebil1} from \cite{O} for $k=1$ and the analog \eqref{kcomdebil1} for $k>1$ from \cite{O2}.

Once again, we only sketch the proof for $k=1$ and refer the reader to the original paper \cite{O} for the details. Let us assume (by homogeneity) that $\|{b}\|_{BMO}=1$. For any $f\in C^\infty_0(\mathbb R^n)$ consider the Calder\'on--Zygmund decomposition of  $f$ at level $\lambda$. We obtain a colection of pairwise disjoint dyadic cubes  $\{Q_j\}=Q_{j}(x_{Q_{j}},r_{j})$ such that
\begin{equation*}\label{CZdec}
\lambda<\frac{1}{|Q_j|}\int_{Q_j}|f(x)|\,dx\leq 2^n\lambda.
\end{equation*}
Define $\Omega=\bigcup_{j}Q_{j}$ and write $f=g+b$, where
$$
g(x)=\left \{ \begin{array}{ll} f(x),&x\in \mathbb R^n\setminus \Omega\\
f_{Q_j},&x\in Q_j,\end{array}\right.
$$
Consider also $h=\displaystyle{\sum_j h_j}$ with $h_j(x)=(f(x)-f_{Q_j})\chi_{Q_j}(x)$ and define, for each  $j$,
$$\widetilde{Q_j}=3Q_j, \quad \widetilde{\Omega}=\displaystyle{\cup_j \widetilde{Q_j}}\quad  \mbox{and} \quad
w_{j}(x)=w(x)\chi_{_{ \mathbb{R}^{n}\setminus 3Q_{j}}}$$
Therefore,
\begin{eqnarray*}
w(\{x\in \mathbb{R}^n:|[b,T]f(x)|>\lambda\})&\leq&w(\{x\in \mathbb{R}^n \setminus \widetilde{\Omega} :|[b,T]g(x)|>\lambda/2\})\\
&+&w(\widetilde{\Omega})\\
&+&w(\{x\in \mathbb{R}^n\setminus
\widetilde{\Omega}:|[b,T]h(x)|>\lambda/2\})\\
&=&I+II+III.
\end{eqnarray*}
Now define  $ \tilde{w}(x):=w(x)\chi_{\mathbb{R}^n\setminus \widetilde{\Omega}}$.  We can deal with the first term by using Chebyshev and the strong result, obtaining that
\begin{equation*}
I\le  \dfrac{c}{\lambda}(pp')^{2p}(r')^{2p-1}\left(\int_{ \mathbb{R}^n \setminus \Omega }|f(x)| \,M_{r}\widetilde{w}(x)\,dx+
\int_{\Omega}|g(x)| \,M_{r}\widetilde{w}(x)\,dx \right).
\end{equation*}
Now, by definition of $\tilde{w}$, away from each $Q_{j}$ the maximal function $M_{r}\tilde{w}$ is almost constant, and hence we obtain that
\begin{equation*}\label{kprimero}
I\leq\frac{c} {\lambda} (pp')^{2p}(r')^{2p-1}\int_{\mathbb{R}^n }|f(x)| \,M_{r}w(x)\,dx
\end{equation*}
The second term is the easy one. It is straightforward to show that
\begin{equation*}
II=w(\widetilde{\Omega})\le \frac{c}{\lambda}\int_{\mathbb R^n} |f(x)|Mw(x)\,dx.
\end{equation*}
Finally, for the third term, we expand the commutator:
\begin{eqnarray*}
III &\leq& w(\{y\in \mathbb{R}^{n} \setminus \widetilde{\Omega}:|\sum_{j}(b(y)-b_{Q_{j}})Th_{j}(y)|>\dfrac{\lambda}{4}\})\\
 &&+w(\{y\in \mathbb{R}^{n}\setminus\widetilde{\Omega}:|\sum_{j}T((b-b_{Q_{j}})h_{j})(y)|>\dfrac{\lambda}{4}\})\\
 &=&A + B.
\end{eqnarray*}
The $A$ term can be bounded by using the properties of the kernel, cancellation of the $h_{j}$'s and standard $BMO$ estimates. We obtain that
\begin{eqnarray*}
A &\leq&  \dfrac{c}{\lambda}\int_{\mathbb{R}^n}|f(x)|\,M_{L(\log L)}w(x)dx\\
 & \le &   cr'\int_{\mathbb{R}^n}\frac{|f(x)|}{\lambda}\,M_rw(x)dx.
\end{eqnarray*}
For the second inequality, we have used that $M_{L(\log L)}w(x)\le c r' M_rw(x)$. Now, for the $B$ term, we can use the known results about weak type boundedness of $T$. It can be proved then (see \cite{O} for details) that
\begin{equation*}
B\leq c\,{(pp')^p(r')^{p-1}}\int_{\mathbb{R}^n}\Phi\left(\dfrac{|f(x)|}{\lambda}\right)\,M_rw(x)dx.
\end{equation*}
for the function $\Phi(t)= t (1+\log^{+}t)$. Combining all estimates together, we obtain the desired result:
\begin{equation*}
w({x\in \mathbb{R}^n: |[b,T]f(x)|>\lambda})\leq
c\,(pp')^{2p}(r')^{2p-1}\int_{\mathbb{R}^n}\Phi\left(\dfrac{|f(x)|}{\lambda}\right)\,M_rw(x)dx.
\end{equation*}
As in the strong case, we refer the reader to \cite{O2} for the details about the analogous result for iterated commutators.

Now we want to optimize inequality \eqref{kcomdebil1} by choosing $p$ and $r$. In the same way as in the strong case, we choose $r=r_{w}=1+\displaystyle{\frac{1}{c_{n}[w]_{A_{\infty}}}}$. This is the best that we can do to control the part with $r'$. In addition, take $p=1+\displaystyle{\frac{1}{\log([w]_{A_{\infty}})}}$. We finish using that $pp'\approx \log([w]_{A_{\infty}})$ and that $A^{1/A}$ is bounded for $A>1$. Define  $\beta= [w]_{A_1}[w]^k_{A_{\infty}}(1+\log^{+}[w]_{A_{\infty}})^{k+1}$ and recall that $\Phi(t)=t(1+\log^{+}t)^k$. Inequality \eqref{kcomdebil1} then becomes
\begin{equation*}
w(\{x\in \mathbb R^n:|T_b^k f(x)|>\lambda\})\leq c_n\,
\beta \,\int_{\mathbb R^n} \Phi
\left(\|{b}\|_{BMO}{\frac{|f(x)|}{\lambda}} \right)\,w(x)dx.
\end{equation*}

\section{Acknowledgement}

The second author is supported by the Spanish Ministry of Science and Innovation grant MTM2009-08934,
the second and third authors are also supported by the Junta de Andaluc\'ia, grant FQM-4745.


\begin{thebibliography}{99}

\bibitem{AP} J.\ Alvarez and C.\ P\'erez, {\it Estimates
with $A_\infty$ weights for various singular integral operators},
Bollettino U.M.I. (7) 8-A (1994), 123--133.


\bibitem{BK1}
R. J. Bagby and D. S. Kurtz, {\it Covering lemmas and the sharp
function,} Proc. Amer. Math. Soc. \textbf{93} (1985), 291-296.


\bibitem{BK2}
R. J. Bagby and D. S. Kurtz, {\it A rearranged good-$\lambda$ 
inequality,} Trans. Amer. Math. Soc. \textbf{293} (1986), 71--81.


\bibitem{B} S. M. Buckley, {\it Estimates for operator norms on weighted spaces
and reverse Jensen inequalities}, Trans. Amer. Math. Soc., {\bf 340}
(1993), no. 1, 253--272.

\bibitem{C}
R.R. Coifman, {\it Distribution function inequalities for singular integrals}, Proc. Nat. Acad. Sci. U.S.A. {\bf 69}
(1972), 2838--2839.

\bibitem{CF}
R.R. Coifman and C. Fefferman, {\it Weighted norm inequalities for
maximal functions and singular integrals}, Studia Math. {\bf 51}
(1974), 241--250.

\bibitem{CMP3} D. Cruz-Uribe, SFO, J.M. Martell and C. P\'erez, {\em
Weights, Extrapolation and the Theory of Rubio de Francia}, Operator Theory: Advances and Applications,
 {\bf 215}, Birkh\"auser/Springer Basel AG, Basel, 2011.


\bibitem{GCRdF} J.\ Garc\'\i a-Cuerva and J.L.\ Rubio de Francia,
{\em Weighted Norm Inequalities and Related Topics}, North Holland
Math. Studies 116, North Holland, Amsterdam, 1985.


\bibitem{Hruscev}
S. Hru\v{s}\v{c}ev, {\it A description of weights satisfying the $A_{\infty}$ condition of Muckenhoupt,} {\it Proc. Amer. Math. Soc.}, {\bf 90}(2), 253--257, 1984.
     	
\bibitem{HP} T. Hyt\"onen and C. Perez, {\em Sharp weighted bounds involving $A_{\infty}$}, preprint.

\bibitem{K}
G. A., Karagulyan,  \emph{Exponential estimates for the {C}alder\'on-{Z}ygmund operator
              and related problems of {F}ourier series}, Mat. Zametki 3,  {\bf 71}, 398--41, 2002.

\bibitem{L} M. Lacey, {\em An $A_p$ --$A_\infty$ inequality for the Hilbert Transform,} preprint. Available at http://arxiv.org/abs/1104.2199

\bibitem{L04}
A.K. Lerner, {\it Weighted rearrangements inequalities for local sharp maximal functions}, Trans. Amer. Math. Soc. (2004), {\bf 357}, (6), 2445 -- 2465.


\bibitem{L2}
A.K. Lerner, {\it Sharp weighted norm inequlities for Littlewood-Paley operators and singular integrals}, Advances in Mathematics (2011), {\bf 226}, (5), 3912 -- 3926.

\bibitem{LOP1}
A. K. Lerner, S. Ombrosi and C. P\'erez, {\it Sharp $A_1$ bounds for
Calder\'on-Zygmund operators and the relationship with a problem of
Muckenhoupt and Wheeden,} International Mathematics Research
Notices, 2008,  \textbf{no. 6}, Art. ID rnm161, 11 pp. 42B20.

\bibitem{LOP3}
A. Lerner, S. Ombrosi and C. P\'erez, {\it $A_1$ bounds for
Calder\'on-Zygmund operators related to a problem of Muckenhoupt and
Wheeden,} Mathematical Research Letters (2009), \textbf{16},
149--156.

\bibitem{LOPTTG}
A. Lerner, S. Ombrosi, C. P\'erez, R. Torres and R. Trujillo-Gonz\'alez, \emph{New maximal functions and multiple weights for the multilinear Calder\'on-Zygmund theory,} Advances in Mathematics (2009), \textbf{220},
1222--1264.

\bibitem{O}
C. Ortiz-Caraballo, {\em Quadratic $A_1$ bounds for commutators of singular
integrals with BMO functions}, to appear in Indiana Univ. Math. J. (2011).

\bibitem{O2}
C. Ortiz-Caraballo, {\em Conmutadores de integrales singulares y pesos $A_1$}, Ph. D. Dissertation, 2011, Universidad de Sevilla.

\bibitem{OPR}
C. Ortiz-Caraballo, C. P\'erez and E. Rela, {\em Local subexponential estimates for classical operators}, Preprint 2011.

\bibitem{Pe1}
C. P\'erez, {\it A course on Singular Integrals and weights}, to appear in Advanced Courses in Mathematics, CRM Barcelona, Birkhauser editors.

\bibitem{P}
C. P\'erez, {\it Weighted norm inequalities for singular integral
operators}, J. London Math. Soc., {\bf 49} (1994), 296--308.

\bibitem{CP1}  C. P\'{e}rez,
{\em Endpoint Estimates for Commutators of Singular Integral
Operators}, Journal of Functional Analysis, {\bf (1) 128} (1995),
163--185.


\bibitem{CP4} C. P\'erez, {\em Sharp estimates for commutators
of singular integrals via iterations of the Hardy-Littlewood maximal
function}, J. Fourier Anal. Appl. 3 (1997), 743-756.


\bibitem{PT} C. P\'erez and R. Trujillo-Gonz\'alez, {\em Sharp weighted estimates for vector-valued singular integrals operators and commutators}, Tohoku Math. J. 55 (2003), 109-129.

\bibitem{Wilson:87}
J. M. Wilson,  {\it Weighted inequalities for the dyadic square function without
dyadic {$A_\infty$}}, Duke Math. J., {\bf 55}(1), 19--50, 1987.

\bibitem{Wilson:89}
J. M. Wilson, {\it Weighted inequalities for the continuous square function}, Trans. Amer. Math. Soc., {\bf 314}(2), 661--692, 1989.

\bibitem{Wilson-LNM}
J. M. Wilson, {\it Weighted Littlewood-Paley theory and exponential-square integrability}, volume 1924 of {\it Lecture Notes in Mathematics.} Springer, Berlin, 2008.



\end{thebibliography}
\end{document}